\documentclass[english, a4paper, 12pt]{article}

\usepackage{babel}
\usepackage{graphicx}
\usepackage{amsmath}
\usepackage{amssymb}
\usepackage{exscale}
\usepackage{hyperref}
\usepackage{times}
\usepackage{subfigure}

\newcommand{\p}{\partial}

\begin{document}

\title{Quasi-Random Discrete Ordinates Method to Radiative Transfer Equation with Linear Anisotropic Scattering}
\author{$^a$Pedro H. A. Konzen\thanks{\texttt{pedro.konzen@ufrgs.br}}, $^b$Leonardo F. Guidi, $^c$Thomas Richter\\
  \small{$^{a,b}$ IME, UFRGS, Porto Alegre, Brazil; $^c$ IAN, OVGU, Magdeburg, Germany}}
\date{}

\maketitle

\section*{Abstract}

The {\it quasi}-random discrete ordinates method (QRDOM) is here proposed for the approximation of transport problems. Its central idea is to explore a {\it quasi} Monte Carlo integration within the classical source iteration technique. It preserves the main characteristics of the discrete ordinates method, but it has the advantage of providing mitigated ray effect solutions. The QRDOM is discussed in details for applications to one-group transport problems with isotropic scattering in rectangular domains. The method is tested against benchmark problems for which DOM solutions are known to suffer from the ray effects. The numerical experiments indicate that the QRDOM provides accurate results and it demands less discrete ordinates per source iteration when compared against the classical DOM.

\noindent\textbf{Keywords:}  radiative transfer equation; anisotropic scattering; {\it quasi}-Monte Carlo method; discrete ordinates method; finite element method.

\section{Introduction}

The modelling of energy transport via radiative transfer is important to many practical high temperature engineering applications \cite{Modest2013a}. To name a few, it is found in the design of industrial furnaces, combustion chambers, or forming processes such as glass and ceramics manufacturing \cite{Frank2004a, Larsen2002a, Viskanta1987a}. Other applications are found in the fields of astrophysics \cite{Meinkohn2002a, Richling2001a}, medical optics \cite{Abdoulaev2003a,  Hielscher1998a, Tarvainen2008a, Wang2007a}, developing of micro-electro-mechanical systems \cite{Knackfuss2006a} and neutron transport \cite{Lewis1984a, Stacey2007a}. 

The steady-state monochromatic radiative transfer equation (RTE, \cite{Modest2013a, Howell2021a}) in homogeneous medium with black boundary conditions is given as follows
\begin{eqnarray}
  &&\forall \pmb{s}\in S^2:~\pmb{s}\cdot\nabla I(\pmb{x},\pmb{s}) + \sigma_tI(\pmb{x},\pmb{s}) \nonumber\\
  &&\qquad\qquad\qquad\qquad= \frac{\sigma_s}{4\pi}\int_{S^2}\Phi(\pmb{s},\pmb{s}')I(\pmb{x},\pmb{s}')\,d\pmb{s}' + kI_b(\pmb{x},\pmb{s}), \forall \pmb{x}\in\mathcal{D},\label{eq:te}\\
  &&\forall \pmb{s}\in S^2, \pmb{n}\cdot\pmb{s} < 0:~I(\pmb{x},\pmb{s}) = I_{\text{in}}(\pmb{x},\pmb{s}), \forall\pmb{x}\in\Gamma,\label{eq:tbc}
\end{eqnarray}
where $S^2 = \{\pmb{s} = (\mu,\eta,\xi)\in\mathbb{R}^3:~\mu^2+\eta^2+\xi^2=1\}$ is the unitary sphere in $\mathbb{R}^3$, centered at origin, $\nabla = \left(\frac{\p}{\p x_1}, \frac{\p}{\p x_2}, \frac{\p}{\p x_3}\right)$ is the gradient operator in $\mathbb{R}^3$, $I(\pmb{x},\pmb{s})$ [$W/(m^2\cdot sr)$] is the radiative intensity at the point $\pmb{x}$ in the domain $\mathcal{D}\subset\mathbb{R}^3$ and in the direction $\pmb{s}\in S^2$, $\sigma_t = k + \sigma_s$ [$m^{-1}$] is the total absorption coefficient, $k$ [$m^{-1}$] and $\sigma_s$ [$m^{-1}$] are, respectively, the homogenized absorption and scattering coefficients, $I_b$ [$W/(m^2\cdot sr)$] and $I_{\text{in}}$ [$W/(m^2\cdot sr)$] are, respectively, the sources in $\mathcal{D}$ and on its boundary $\Gamma$, $\pmb{n}$ is the unit outer normal on $\Gamma$. The scattering phase function is assumed to be given as
\begin{eqnarray}\label{eq:Phi}
  &&\Phi(\pmb{s},\pmb{s}') = a_0 + a_1\pmb{s}\cdot\pmb{s}',
\end{eqnarray}
and anisotropic scattering is stated by setting $a_1\neq 0$.

The discrete ordinates method (DOM) is one of the most widely used techniques to solve Eqs.~\ref{eq:te}-\ref{eq:tbc} (see, for instance, \cite[Ch. 17]{Modest2013a}). It consists in approximating the integral term in the right-hand side of Eq.~\ref{eq:te} by using an appropriate quadrature set $\{\pmb{s}_i, \omega_i\}_{i=1}^M$. This leads to the approximation of the integro-differential equation by a system of partial differential equations on the discrete ordinates $\{I(\pmb{x},\pmb{s}_i)\}_{i=1}^M$, which can be solved by a variety of classical discretization methods. Unfortunately, for transport problems with discontinuities in $I_b$, $\sigma_{t}$, $\sigma_{s}$  or on non-convex geometries, the DOM approximation may produce unrealistic oscillatory solutions known as the ray effects \cite{Chai1993a, Morel2003a}.

At the expense of additional computational costs, ray effects can be mitigated by increasing the number of discrete ordinates \cite{Li2003a}. It is also well known that the quality of the DOM solution depends on the choice of the quadrature set and is problem depended \cite{AbuShumays2001a,Barichello2016a,Hunter2013a,Koch2004a}. Therefore, Adaptive Discrete Ordinates schemes \cite{Jarrel2010a, Stone2007a} have been proposed. Alternatively, the Frame Rotation Method \cite{Tencer2016a} was developed to not depend on the choice of the quadrature rule and provide solutions invariant under arbitrary rotations of the reference frame. In short, this method computes the solution as the average of DOM solutions with random rotations of a given quadrature rule. One of its disadvantages is that the frame rotations do not preserve the symmetries of the discrete directions related to the domain boundary, forcing approximations of the solution on the boundary that may decrease the overall accuracy.

The {\it Quasi}-Random Discrete Ordinates Method (QRDOM, \cite{Konzen2019a}) is an alternative technique that maintains most characteristics of the DOM with mitigated ray effects. Its central idea is to explore a {\it quasi}-Monte Carlo integration \cite{Leobacher2014a} within the classical source iteration technique. Firstly developed for the case of isotropic scattering, it is recently being extended and tested for anisotropic problems \cite{Konzen2022a, Konzen2023a}. The main objective of this paper is to discuss and test the extension of the QRDOM to solve the transport problem Eqs.~\ref{eq:te}-\ref{eq:tbc} in a rectangular domain $\mathcal{D}\subset\mathbb{R}^2$ with linear anisotropic scattering modeled by the phase function given in Eq.~\ref{eq:Phi}. Verification of the proposed extension is performed by applications to problems with manufactured solutions.

\section{QRDOM for Linear Anisotropic Scattering}

\noindent The QRDOM to transport problems with isotropic scattering is explained in detail in \cite{Konzen2019a}. Here, it is extended to solve the transport problem Eqs.~\ref{eq:te}-\ref{eq:tbc} with linear anisotropic scattering phase function, Eq.~\ref{eq:Phi}. With this in mind, the integral term in Eq.~\ref{eq:te} can be rewritten as follows
\begin{eqnarray}
  \Psi(\pmb{x},\pmb{s})&:=&\frac{1}{4\pi}\int_{S^2}\Phi(\pmb{s},\pmb{s}')I(\pmb{x},\pmb{s}')d\pmb{s}' \label{eq:Psi}\\
  &=&\frac{a_0}{4\pi}\int_{S^{2}}I(\pmb{x},\pmb{s}')d\pmb{s}' + a_1\sum_{i=1}^3\frac{s_i}{4\pi}\int_{S^2}s_i'I(\pmb{x},\pmb{s}')d\pmb{s}'\\
  &=:&a_0\Psi_0(\pmb{x}) + a_1\sum_{i=1}^3s_i\Psi_i(\pmb{x}),
\end{eqnarray}
where $\Psi_0(\pmb{x})$ [$W/m^2$] is the radiation density, and $\Psi_i(\pmb{x})$ [$W/m^2$] the $i$-th partial current density.

The computational domain for the classical source iteration approximation with Eqs.~\ref{eq:te}-\ref{eq:tbc} is assumed to be rectangular $\mathcal{D} = [a, b]\times [c, d]\in \mathbb{R}^2$. The iteration reads
\begin{eqnarray}
  &&\forall \pmb{s}\in S^2:~s_1\frac{\p}{\p x_1}I^{(l)} + s_2\frac{\p}{\p x_2}I^{(l)} + \sigma_tI(\pmb{x},\pmb{s}) = \sigma_s\Psi^{(l-1)}(\pmb{x},\pmb{s}) + kI_b(\pmb{x})\\
  &&\forall \pmb{s}\in S^2, \pmb{n}\cdot\pmb{s} < 0:~I^{(l)}(\pmb{x},\pmb{s}) = I_{\text{in}}(\pmb{x},\pmb{s})
\end{eqnarray}
where, $l = 1, 2, 3, \ldots$, $\Psi^{(0)}$ is the initial approximation of Eq.~\ref{eq:Psi} and $\Psi^{(l)}(\pmb{x},\pmb{s}) = a_0\Psi_0^{(l)}(\pmb{x}) + a_1s_1\Psi_1^{(l)}(\pmb{x}) + a_1s_2\Psi_2^{(l)}(\pmb{x})$.  A iteration stop criteria based on $\Psi^{(l)}$ and $\Psi^{(l-1)}$ is considered. 

\begin{figure}[h]
  \centering
  \includegraphics[width=0.6\textwidth]{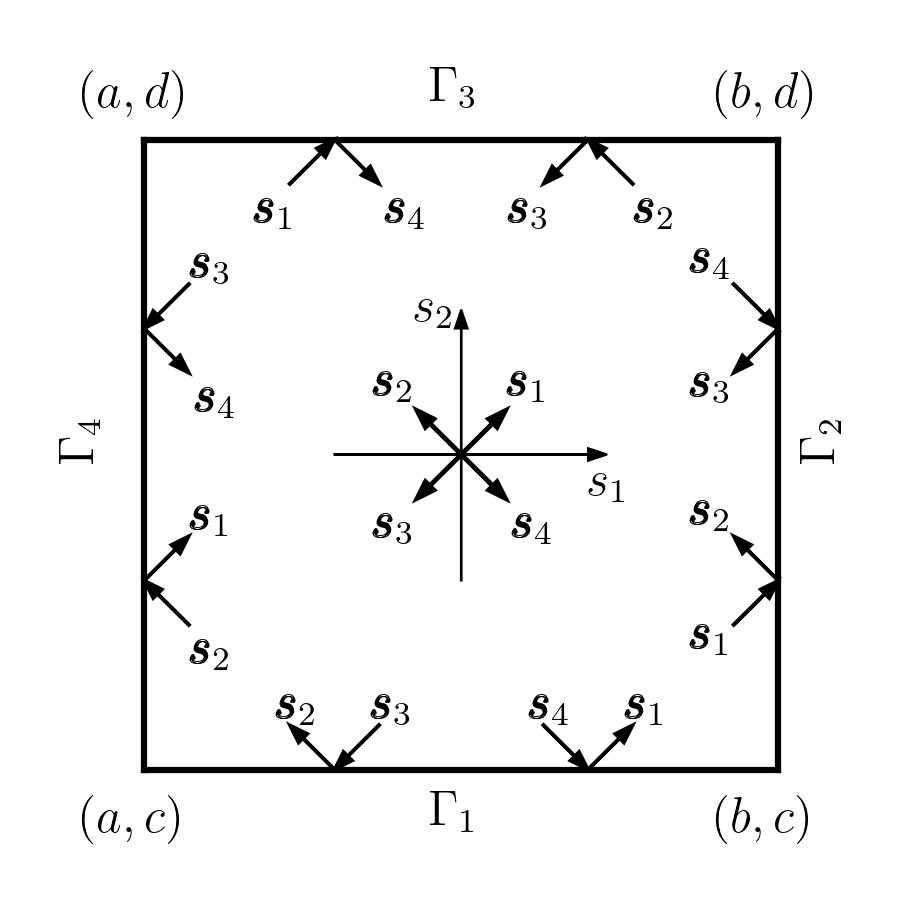}
  \caption{Illustration of the rectangular domain with reflected boundary directions.}
  \label{fig:rect_domain}
\end{figure}

The main idea of the QRDOM is to use a {\it quasi}-Monte Carlo integration to compute $\Psi^{(l)}$. A {\it quasi}-random discrete ordinates generator $qr_1:\mathbb{N}\to S^2$ is built from the bi-dimensional $\mathcal{S}_{2,3}$ {\it quasi}-random reverse Halton sequence number generator $rh:\mathbb{N}\to (0,1)^2$ \cite{Vandewoestyne2006a}. More explicitly, $qr_1(i) = \pmb{s}_{i,1}$ is given as
\begin{eqnarray}
  s_{i,1,1}&=&\cos(\phi_i)\sin(\theta_i)\\
  s_{i,1,2}&=&\sin(\phi_i)\sin(\theta_i)\\
  s_{i,1,3}&=&\cos(\theta_i)
\end{eqnarray}
where, the $\theta_i = \arccos\left(1-rh_1(i)\right)$ is the polar angle and $\phi_i = rh_2(i)\pi/2$ is the azimuth angle, $i = 1, 2, 3, \ldots$. The $qr_1$ generates {\it quasi}-random directions on the first octant of $S^2$, and the boundary reflected directions (Fig.~\ref{fig:rect_domain}) are given by $\pmb{s}_{i,2} = (-s_{i,1,1}, s_{i,1,2}, s_{i,1,3})$, $\pmb{s}_{i,3} = (-s_{i,1,1}, -s_{i,1,2}, s_{i,1,3})$, $\pmb{s}_{i,4} = (s_{i,1,1}, -s_{i,1,2}, s_{i,1,3})$. See \cite[Table 1 and Fig. 2]{Konzen2019a} for more details about the generated directions.

The QRDOM iteration then reads
\begin{eqnarray}
  j=1,2,3,4:&& s_{i,j,1}\frac{\p I_{i,j}^{(l)}}{\p x} + s_{i,j,2}\frac{\p I_{i,j}^{(l)}}{\p y} + \sigma_tI_{i,j}^{(l)} = \sigma_s\Psi_{i,j}^{(l-1)} + kI_b(\pmb{x}), \forall \pmb{x}\in\mathcal{D}, \label{eq:qrtp_rect_eq}\\
  && I_{i,1}^{(l)} = I_{i,2}^{(l)} = I_{\text{in}},~ \text{on}~\Gamma_1:=[a,b]\!\times\!\{0\}, \label{eq:qrtp_rect_bc1}\\ 
  && I_{i,2}^{(l)} = I_{i,3}^{(l)} = I_{\text{in}},~ \text{on}~\Gamma_2:=\{b\}\!\times\![c, d],  \label{eq:qrtp_rect_bc2}\\
  && I_{i,3}^{(l)} = I_{i,4}^{(l)} = I_{\text{in}}, ~ \text{on}~\Gamma_3:=[a,b]\!\times\!\{d\},  \label{eq:qrtp_rect_bc3}\\
  && I_{i,1}^{(l)} = I_{i,4}^{(l)} = I_{\text{in}}, ~ \text{on}~\Gamma_4:=\{a\}\!\times\![c,d], \label{eq:qrtp_rect_bc4}
\end{eqnarray}
for each $i = m^{(l-1)}, m^{(l-1)}+1, m^{(l-1)}+2, \dotsc, m^{(l)}-1$, where
\begin{eqnarray}
  &&\Psi_{i,j}^{(l)} = a_0\Psi_0^{(l)} + a_1s_{i,j,1}\Psi_1^{(l)} + a_1s_{i,j,2}\Psi_2^{(l)}
\end{eqnarray}
and, $m^{(l)}-1$ is the lower index of the quasi-random sequence greater than $m^{(l-1)}$ for which $\Psi^{(l)}$ converges to a given precision. The approximation of the integral terms for the next iteration is accumulatively computed as
\begin{eqnarray}
  &&\Psi_0^{(l)} = \frac{a_0}{2M^{(l)}}\sum_{k=m^{(l-1)}}^{m^{(l)}-1}\sum_{j=1}^4I_{k,j}^{(l)}\label{eq:psi_0}\\
  &&\Psi_1^{(l)} = \frac{a_1}{2M^{(l)}}\sum_{k=m^{(l-1)}}^{m^{(l)}-1}\sum_{j=1}^4s_{k,j,1}I_{k,j}^{(l)}\\
  &&\Psi_2^{(l)} = \frac{a_1}{2M^{(l)}}\sum_{k=m^{(l-1)}}^{m^{(l)}-1}\sum_{j=1}^4s_{k,j,2}I_{k,j}^{(l)}\label{eq:psi_2}
\end{eqnarray}
where $M^{(l)} := m^{(l)}-m^{(l-1)}$ is the number of {\it quasi}-random ordinates used in the $l$-th iteration.

\textbf{Spatial Discretization.} The QRDOM requires the application of a numerical method to solve the iterate problem Eqs.~\ref{eq:qrtp_rect_eq}-\ref{eq:qrtp_rect_bc4}. Here, the standard finite element method (FEM) with the streamline upwind Petrov-Galerkin stabilization (SUPG) is applied (see, for instance, \cite{Kanschat1998a} for a similar application). We define the spaces $V_{i,j} = \{v\in L^2(\mathcal{D}:~ \pmb{s}_{i,j}\cdot\nabla v \in L^2(\mathcal{D}))\}$. The related weak problem is stated with the SUPG trial function $v_{i,j,\delta} := v + \delta \pmb{s}_{i,j}\nabla v$, with $v\in V_{i,j}$ and with the stabilization parameter 
\begin{eqnarray}
  &&\delta = \left(\frac{c_1^2}{h^2} + \sigma_t\right)^{-\frac{1}{2}},
\end{eqnarray}
where $c_1=2$. From Eq.~\ref{eq:qrtp_rect_eq} it follows 
\begin{eqnarray}
  &&\left(\pmb{s}_{i,j}\cdot\nabla I_{i,j}^{(l)}, v+\delta\pmb{s}_{i,j}\nabla v\right)_{\mathcal{D}} + \left(\sigma_t I_{i,j}^{(l)}, v_{i,j,\delta}\right)_{\mathcal{D}} \nonumber\\
  &&\qquad\qquad\qquad\qquad= \left(\sigma_s\Psi^{(l-1)} + \kappa I_b, v_{i,j,\delta}\right)_{\mathcal{D}}, \forall v\in V_{i,j},
\end{eqnarray}
where $\left(\cdot,\cdot\right)_{\mathcal{D}}$ denotes the standard $L^2$ inner product. By applying the Green theorem, it follows
\begin{eqnarray}
  &&-\left(I_{i,j}^{(l)}, \pmb{s}_{i,j}\cdot\nabla v\right)_{\mathcal{D}} + \left(I_{i,j}^{(l)}, \pmb{s}_{i,j}\cdot \pmb{n} v\right)_{\Gamma_{i,j}^+} + \left(\pmb{s}_{i,j}\cdot\nabla I_{i,j}^{(l)}, \delta\pmb{s}_{i,j}\nabla v\right)_{\mathcal{D}} \nonumber\\
  &&\qquad\qquad\qquad+ \left(\sigma_t I_{i,j}^{(l)}, v_{i,j,\delta}\right)_{\mathcal{D}} = \left(\sigma_s\Psi^{(l-1)} + \kappa I_b, v_{i,j,\delta}\right)_{\mathcal{D}} \nonumber\\
  &&\qquad\qquad\qquad- \left(I_{\text{in}}, \pmb{s}_{i,j}\cdot \pmb{n} v\right)_{\Gamma_{i,j}^-}, \forall v\in V_{i,j},\label{eq:qr_weak}
\end{eqnarray}
where $\Gamma_{i,j}^+ := \{x\in\p\mathcal{D}:~\pmb{s}_{i,j}\cdot n > 0\}$ and $\Gamma_{i,j}^- := \{x\in\p\mathcal{D}:~\pmb{s}_{i,j}\cdot n < 0\}$. For the sake of simplicity, it will be assumed that $a(I_{i,j}, v_{i,j,\delta})$ denotes the left-hand-side and $l(v_{i,j,\delta})$ the right-hand-side of Eq.~\ref{eq:qr_weak}.

From the weak form Eq.~\ref{eq:qr_weak}, the discrete finite element problem reads: for each $i$ find $I_{i,j,h}\in V_{i,j,h}$, such that
\begin{eqnarray}\label{eq:qr_fem}
  &&a\left(I_{i,j,h}, v_{i,j,\delta,h}\right) = l(v_{i,j,\delta,h}),\quad\forall v\in V_{i,j,h},
\end{eqnarray}
for each $j=1,2,3,4$, where $V_{i,j,h}\subset V_{i,j}$ is the finite element space of quadratic elements $Q_1$ built on a regular mesh. 

\textbf{Implementation Details.} Usually the QRDOM will require Eq.~\ref{eq:qr_fem} to be solved hundreds of times at each source iteration, i.e. $M^{(l)}$ is expected to be of order of $10^2$ or $10^3$ to reach a commonly accepted tolerance. These problems can be solved in parallel, but a race condition occurs in the computation of $\Psi^{(l)}$. Care must be exercised during evaluation of $\Psi_0^{(l)}$, $\Psi_1^{(l)}$ and $\Psi_2^{(l)}$ since its terms may be acquired by parallelization's procedure in an order different from that of the {\it quasi}-random sequence as defined in Eqs.~\ref{eq:psi_0}-\ref{eq:psi_2}.

In the next section, numerical results for test cases are presented. They have been achieved by implementing the QRDOM with the help of the finite element toolkit Gascoigne 3D \cite{Becker2021a}, following a parallel implementation with message passing interface (MPI) and multiprocessing (MP). The linear system Eq.~\ref{eq:qr_fem} has been solved with the GMRES method with the relative unpreconditioned residual tolerance set to $10^{-12}$.

\section{Results}

\noindent Here the QRDOM is applied to two manufactured solutions of radiative transfer problems with linear anisotropic scattering. In the first problem the radiative intensity $I$ does not depend on the direction $\pmb{s}$. In the second problem, this dependence is considered and QRDOM solutions are verified for a range of $\sigma_s\in [0.1, 5]$.

\begin{table}[!htbp] 
\caption{Selected functionals on the radiation density computed for Problem 1 with the QRDOM for several mesh global refinements.}
\vspace{12pt}
\centering{}
\renewcommand{\arraystretch}{1.2}
\begin{tabular*}{\textwidth}{@{\extracolsep{\fill}}c|cc|ccc}\hline
\#cells & $\epsilon$ & $F(\Psi_0)$ & $\Psi_0(0.1,0.1)$ & $\Psi_0(0.5,0.5)$ & $\Psi_0(0.9,0.9)$\\\hline
$16\times 16$ & 5.1e-03 & 0.999994 & 1.54049 & 0.99999 & 1.54049\\\hline 
$32\times 32$ & 6.1e-04 & 0.999994 & 1.55876 & 0.99999 & 1.55876\\\hline 
$64\times 64$ & 7.5e-05 & 0.999994 & 1.55797 & 0.99999 & 1.55797\\\hline 
$128\times 128$ & 1.1e-05 & 0.999994 & 1.55904 & 0.99999 & 1.55904\\\hline 
Exact & -x- & 1.000000 & 1.55902 & 1.00000 & 1.55902\\\hline 
\end{tabular*}\label{tab:selectedFuns}
\end{table}

\begin{figure}[!htbp] 
\vspace{-2pt}
\begin{center}
\includegraphics[width=0.7\textwidth]{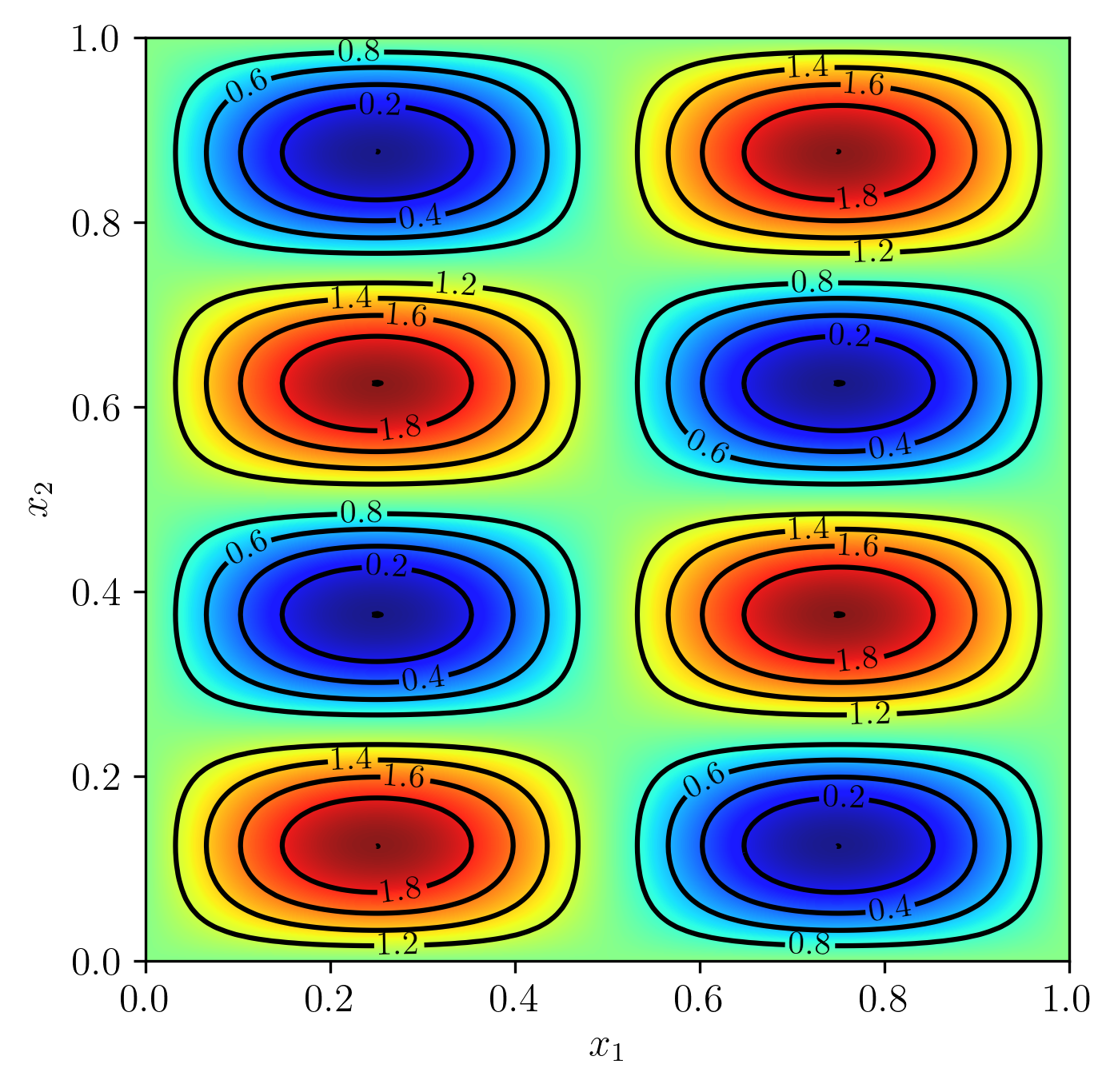}
\caption{QRDOM radiation density $\Psi_0$ for Problem 1 computed with a mesh of $128\times 128$ cells.}
\label{fig:contourf}
\end{center}
\end{figure}

\begin{figure}[!htbp] 
\vspace{-2pt}
\begin{center}
\includegraphics[width=0.7\textwidth]{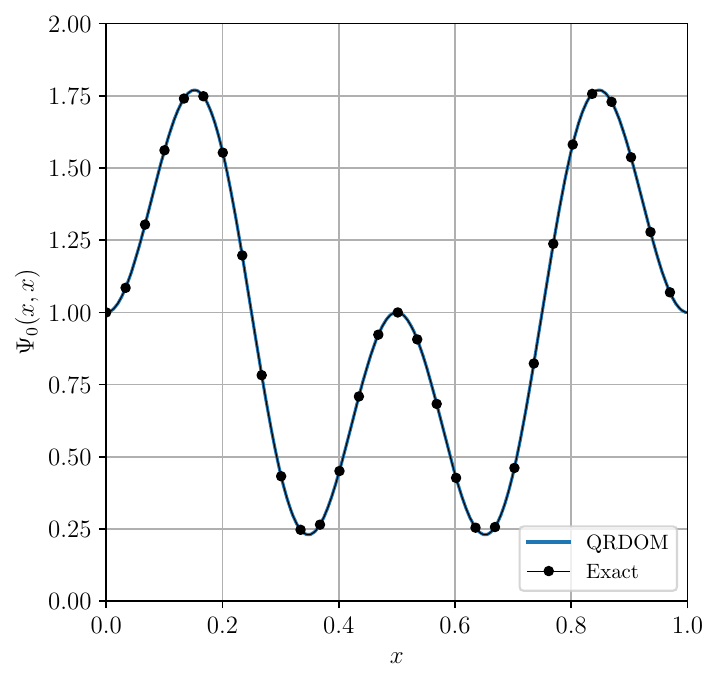}%
\caption{QRDOM radiation density $\Psi_0(x_1,x_2)$ for Problem 1 computed on the identity line $x_2=x_1$.}
\label{fig:tline}%
\end{center}
\end{figure}

\textbf{Problem 1.} Following the method of manufactured solutions, we assume the exact radiative intensity
\begin{eqnarray}\label{eq:exact_sol}
  &&\hat{I}(x_1,x_2) = 1 + \sin(2\kappa\pi x_1)\sin(2\sigma_s\pi x_2),
\end{eqnarray}
with $\sigma_s = 2\kappa$ and $\kappa = 1$, as proposed by \cite{Badri2018a}. The phase function coefficients are $a_0 = 1$ and $a_1 = 1/2$, compare Eq.~\ref{eq:Phi}. Substituting $\hat{I}(x_1,x_2)$ in problem Eqs.~\ref{eq:te}-\ref{eq:tbc} and assuming a rectangular domain $\mathcal{D} = [0, 1]\times [0, 1]$ (2D symmetry), the source is found to be
\begin{eqnarray}
  kI_b(x_1,x_2) &=& 2\pi\kappa s_{1}\cos(2\pi\kappa x_1)\sin(2\pi\sigma_s x_2) \nonumber\\
  &+& 2\pi\sigma_s s_2\sin(2\pi\kappa x_1)\cos(2\pi\sigma_s x_2) + \kappa \hat{I}(x_1, x_2).
\end{eqnarray}
The exact radiation density is
\begin{eqnarray}
  &&\hat{\Psi}_0(x_1,x_2) = 1 + \sin(2\pi\kappa x_1)\sin(2\pi\sigma_s x_2),
\end{eqnarray}
and the solution has null partial current densities.

The stopping criterium for the inner QRDOM iterations and also for the source iterations is set to
\begin{eqnarray}\label{eq:stop_criteria}
  &&\frac{\left|F\left(\Psi_0^{(l)}\right) - F\left(\Psi_0^{(l-1)}\right)\right|}{\left|F\left(\Psi_0^{(l)}\right)\right|} < 10^{-5},
\end{eqnarray}
where, $F$ denotes the target functional
\begin{eqnarray}\label{eq:tf}
 &&F(\Psi) = \frac{1}{|\mathcal{D}|}\int_\mathcal{D}\Psi(\pmb{x})\,d\pmb{x}.
\end{eqnarray}

Table \ref{tab:selectedFuns} contains selected functionals on the radiation density computed with the QRDOM for several global mesh refinements. The error $\epsilon$ is the $L^2$ norm of the difference between the exact and the QRDOM computed radiation density $\Psi_0$. It is notable that the results have good precision and are in accordance to the tolerance of $10^{-5}$ set as stop criteria. Fig.~\ref{fig:contourf} shows the QRDOM computed radiation density with a uniform mesh of $128\times 128$ cells. For this problem, the solution has null partial current densities and the QRDOM provides them in the machine precision of about $10^{-16}$.

A further comparison between the exact and the QRDOM computed radiation density is shown in Fig.~\ref{fig:tline}. The blue solid line is the QRDOM solution on the identity line $x_2=x_1$, and the black dotted line the exact solution. Following the good precision of the QRDOM already observed in Table \ref{tab:selectedFuns}, it is also observable that the symmetry of the solution is well preserved by the QRDOM.

\textbf{Problem 2.} The manufactured solution now considered is a function of both the space and the direction, as in \cite{LeHardy2016a}. The radiative intensity is assumed to be
\begin{eqnarray}
  &&\hat{I}(\pmb{x},\pmb{s}) = (1+s_1)e^{-\kappa x_1 - \sigma_s x_2},\quad \mathcal{D} = [0,1]\times [0,1].
\end{eqnarray}
By substituting it in Eqs.~\ref{eq:te}-\ref{eq:tbc}, the source is found to be
\begin{eqnarray}
  &&\kappa I_b(\pmb{x},\pmb{s}) = \left[(\kappa - \kappa s_1 - \sigma_s s_2)(1 + s_1)+\frac{5}{6}\sigma_s s_1\right]e^{-\kappa x_1 - \sigma_s x_2}.
\end{eqnarray}
The exact radiation density is
\begin{eqnarray}
  &&\hat{\Psi}_0(x_1,x_2) = e^{-\kappa x_1 - \sigma_s x_2},
\end{eqnarray}
and the partial current densities are $\hat{\Psi}_1(\pmb{x}) = \frac{1}{3}\hat{\Psi}_0(\pmb{x})$ and $\hat{\Psi}_2(\pmb{x}) = 0$. As before, the stopping criterium for the inner QRDOM iterations and also for the source iterations is set by Eq.~\ref{eq:stop_criteria}.

\begin{table}[!hbtp] 
\caption{Comparison between QRDOM ($\Psi$, mesh with $128\times 128$ cells) and exact ($\hat{\Psi}$) solutions for Problem 2 with $\kappa=0.1$.}
\vspace{12pt}
\centering{}
\renewcommand{\arraystretch}{1.2}
\begin{tabular*}{\textwidth}{@{\extracolsep{\fill}}c|c|cc|cc|c}\hline
  $\sigma_s$ & $\epsilon$ & $F(\Psi_0)$ & $F(\hat{\Psi}_0)$ & $F(\Psi_1)$ & $F(\hat{\Psi}_1)$ & $F(\Psi_2)$\\\hline
  0.1 & 3.8e-06 & 0.905592 & 0.905592 & 0.299835 & 0.301864 & 5.7e-10 \\\hline
  0.9 & 4.2e-05 & 0.627471 & 0.627471 & 0.209219 & 0.209157 & -5.1e-08\\\hline
  2.5 & 9.6e-06 & 0.349412 & 0.349405 & 0.116567 & 0.116468 & 1.3e-08\\\hline
  5.0 & 2.1e-05 & 0.189055 & 0.189043 & 0.062626 & 0.063014 & 4.8e-08\\\hline
\end{tabular*}\label{tab:prob2_tf}
\end{table}

\begin{figure}[!htbp] 
\vspace{-2pt}
\begin{center}
  \includegraphics[width=0.49\textwidth]{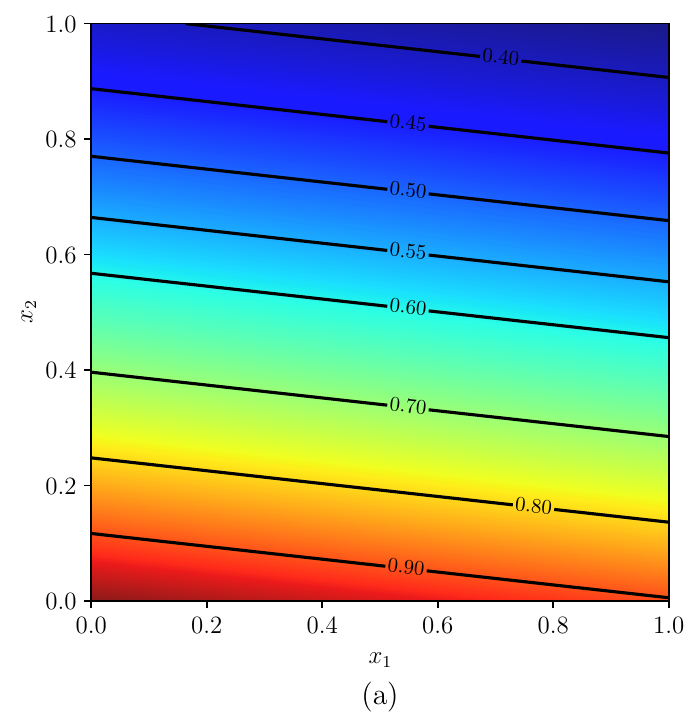}~%
  \includegraphics[width=0.49\textwidth]{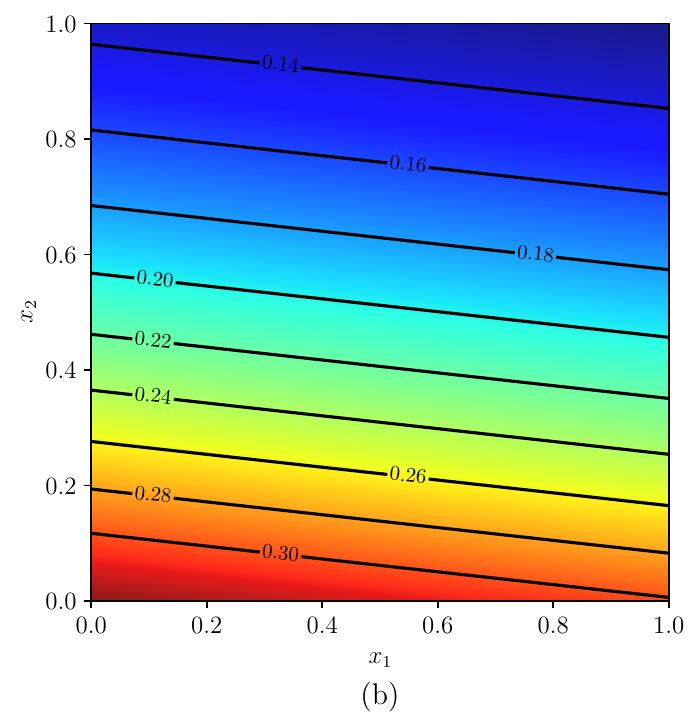}
\caption{QRDOM solution ($128\times 128$ mesh) for Problem 2 with $\kappa=0.1$ and $\sigma_s=0.9$. (a) $\Psi_0$. (b) $\Psi_1$.}
\label{fig:prob2_contourf}%
\end{center}
\end{figure}

\begin{figure}[!htbp] 
\vspace{-2pt}
\begin{center}
\includegraphics[width=0.7\textwidth]{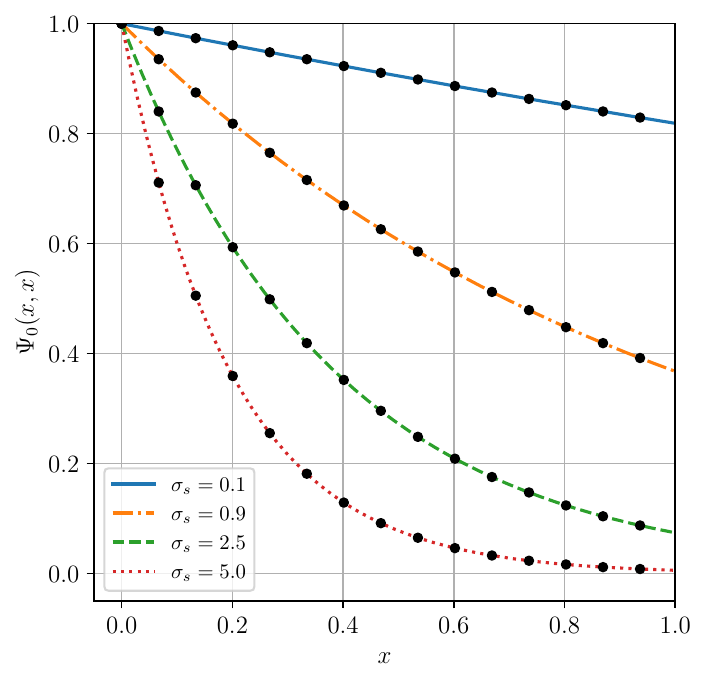}%
\caption{QRDOM radiation density $\Psi_0(x_1,x_2)$ solutions (lines) on the identity line $x_1=x_2$ {\it versus} the manufactured solution for Problem 2 with different values of $\sigma_s$.}
\label{fig:prob2_idline}%
\end{center}
\end{figure}

Table \ref{tab:prob2_tf} presents the values of the target functional Eq.~\ref{eq:tf} of the QRDOM and the exact solutions for Problem 2 with $\kappa=0.1$. The QRDOM solutions were computed with a uniform mesh of $128\times 128$ cells and, as before, $\epsilon$ denotes the $L^2$ error computed on the base of the radiation density. One can observe that the values of the target functional of the QRDOM computed radiation densities ($F(\Psi_0)$) have a good accuracy in the order of the stop criteria tolerance. This is expected once $F(\Psi_0)$ is used by the method to track the convergence. A drop in accuracy can be observed for the computed partial current densities, but it is equivalent to the expected from the classical DOM method with similar parameters ($\sim 100-500$ ordinates directions on the first octant of $S^2$). Not shown in this table, is the exact value $F(\hat{\Psi}_2) = 0$.

Fig.~\ref{fig:prob2_contourf} shows the QRDOM solution with a $128\times 128$ mesh for the Problem 2 with $\kappa=0.1$ and $\sigma_s=0.9$. Fig.~\ref{fig:prob2_contourf}(a) shows the computed radiation density $\Psi_0$ and  Fig.~\ref{fig:prob2_contourf}(b) the partial current density $\Psi_1$. As it is expected, the isolines of the QRDOM solution are straight lines.

Fig.~\ref{fig:prob2_idline} shows the comparison between the QRDOM radiation density $\Psi_0$ solution and the exact solution of Problem 2. The QRDOM solutions were obtained with a mesh of $128\times 128$ cells and the total absorption coefficient set to $\kappa = 0.1$. For all considered values of the scattering coefficient $\sigma_s= 0.1, 0.9, 2.5, 5.0$, the $L^2$ error between the method and the exact solutions are $\epsilon \approx 10^{-5}$, in accordance with the assumed stop criteria.

\section{Summary}

\noindent This paper discussed the extension and application of the QRDOM, originally proposed for transport problems with isotropic scattering, to radiative transfer problems with linear anisotropic scattering. Instead of just the radiation density, here the application requires the {\it quasi}-Monte Carlo integration of the partial current densities. Other aspects of the method's implementation are preserved, allowing the use of a similar parallel MPI/MP strategy as in the original version.

In order to validate the novel developments, the method was applied to manufactured solutions of problems stated in a rectangular domain with 2D symmetry and black boundaries. The good results achieved by the method indicate its potential for radiative transfer problems with anisotropic scattering.

Further developments include testing the QRDOM for similar problems with discontinuities in the radiation source, absorption or scattering coefficients. These are specially important test cases to evaluate its potential mitigation of the ray effects, an issue known to be observed with the classical DOM. Other developments may accelerate the QRDOM convergence by exploring the source iteration and the properties of the {\it quasi}-random sequence generator.

\bibliographystyle{plain}

\end{document}